\documentclass{amsart}

\usepackage{amssymb}
\usepackage{amstext}
\usepackage{amsmath}
\usepackage{amscd}
\usepackage{amsthm}
\usepackage{amsfonts}
\usepackage{enumerate}
\usepackage{graphicx}
\usepackage{color}
\usepackage{here} % 図や表を強制的に出力させる
\usepackage{bm} % 太字ベクトル
\usepackage{mathrsfs}
\usepackage{mathtools}
\usepackage{latexsym}
\usepackage{booktabs}
\usepackage{fancyhdr}
\usepackage{subcaption}
\usepackage{tikz}
\usepackage{tikz-cd}
\usetikzlibrary{arrows}
\usetikzlibrary{decorations.markings}
\usetikzlibrary{patterns,decorations.pathreplacing}
\usepackage{pdfpages}
\usepackage{hyperref}
\hypersetup{
  colorlinks=true,
  citecolor=red,
  linkcolor=blue
}
\theoremstyle{plain}
\newtheorem{theorem}{Theorem}[section]
\newtheorem{lemma}[theorem]{Lemma}
\newtheorem{proposition}[theorem]{Proposition}
\newtheorem{corollary}[theorem]{Corollary}
\newtheorem{conjecture}[theorem]{Conjecture}

\theoremstyle{definition}

\newtheorem{example}[theorem]{Example}

\theoremstyle{remark}

\newcommand{\intposit}{\mathbb{Z}_{> 0}}

\newcommand{\calF}{\mathcal{F}}

\newcommand{\calH}{\mathcal{H}}

\newcommand{\calZ}{\mathcal{Z}}

\newcommand{\ZZ}{\mathbb{Z}}

\newcommand{\RR}{\mathbb{R}}

\newcommand{\kk}{\Bbbk}

\newcommand{\ab}{\mathbf{a}}

\newcommand{\bb}{\mathbf{b}}

\newcommand{\xb}{\mathbf{x}}
\newcommand{\yb}{\mathbf{y}}

\newcommand{\zb}{\mathbf{z}}

\DeclareMathOperator{\ini}{in}

\DeclareMathOperator{\Span}{span}

\title{Graded Ehrhart theory for type~$A$ alcoved polytopes}
\date{}

\author[M. Konoike]{Masato Konoike}
\address{Department of Pure and Applied Mathematics, Graduate School of Information
Science and Technology, The University of Osaka, Osaka, Japan}
\email{kounoike-m@ist.osaka-u.ac.jp}

\author[S. Miyashita]{Sora Miyashita}
\address{Department of Pure and Applied Mathematics, Graduate School of Information
Science and Technology, The University of Osaka, Osaka, Japan}
\email{u804642k@ecs.osaka-u.ac.jp}

\author[K. Tani]{Koichiro Tani}
\address{Department of Pure and Applied Mathematics, Graduate School of Information
Science and Technology, The University of Osaka, Osaka, Japan}
\email{tani-k@ist.osaka-u.ac.jp}

\subjclass[2020]{Primary 52B20; Secondary 05E40, 13F65, 13H10}
\keywords{Graded Ehrhart theory, alcoved polytope, harmonic algebra, coordinate compression, initial algebra, canonical ideal}

\begin{document}
\begin{abstract}
    In this paper, we study harmonic algebras of type~$A$ alcoved polytopes.
    For each type~$A$ alcoved polytope, we compute the initial algebra of its harmonic algebra with respect to every monomial order and prove that every such initial algebra is a normal and Cohen--Macaulay affine semigroup ring.
    In particular, the harmonic algebra is finitely generated, normal, and Cohen--Macaulay.
    When the polytope is full-dimensional, we construct a canonical ideal whose initial ideal is the canonical ideal of the initial algebra for every monomial order.
\end{abstract}

\maketitle

\section{Introduction}\label{sec:introduction}
Ehrhart theory studies the lattice points in the dilates of a lattice polytope. For a lattice polytope $P\subseteq\RR^n$, its \emph{Ehrhart polynomial} and \emph{Ehrhart series} are
\[
 i_P(m)=|mP\cap\ZZ^n|,
 \qquad
 \mathrm E_P(t)=\sum_{m\geq0}i_P(m)t^m.
\]
The rationality of this series and Ehrhart--Macdonald reciprocity have
algebraic interpretations through the Ehrhart ring of $P$: finite
generation gives rationality, while Cohen--Macaulayness and the
description of the canonical module by interior lattice points give
reciprocity. Graded Ehrhart theory, introduced by Reiner and
Rhoades~\cite{ReinerRhoades}, asks for analogous statements for a
$q$-refinement of lattice point enumeration constructed using orbit
harmonics.

Throughout this paper, let $\kk$ be a field of characteristic zero. To a finite lattice set $\calZ\subseteq\ZZ^n$, Reiner and Rhoades associate a harmonic space $V_{\calZ}\subseteq\kk[\yb]$, where
$\yb=(y_1,\ldots,y_n)$, satisfying $\dim_{\kk}V_{\calZ}=|\calZ|$. Its Hilbert series agrees with that of the point-orbit ring $R(\calZ)$; we recall both constructions in Section~\ref{section:basics of harmonic algebras}. Setting $\calZ_m=mP\cap\ZZ^n$, define the \emph{graded Ehrhart function} and the \emph{graded Ehrhart series} of $P$ by
\[
 i_P(m;q)=\operatorname{Hilb}(V_{\calZ_m},q),
 \qquad
 \mathrm E_P(t,q)=\sum_{m\geq0}i_P(m;q)t^m.
\]
The coefficients of $i_P(m;q)$ are nonnegative integers and
$i_P(m;1)=i_P(m)$. Thus $\mathrm E_P(t,q)$ specializes at $q=1$ to the
ordinary Ehrhart series.
Reiner and Rhoades assemble these harmonic spaces into the \emph{harmonic algebra}
\[
 \calH_P=\bigoplus_{m\geq0}y_0^mV_{\calZ_m}\subseteq\kk[y_0,\yb],
 \qquad
 \deg y_0=(1,0),\quad \deg y_i=(0,1),
\]
where $y_i$ ranges over the variables in $\yb$.
Its bigraded Hilbert series is $\mathrm E_P(t,q)$. For $m\geq1$, set $\calZ_m^\circ=\operatorname{relint}(mP)\cap\ZZ^n$, where $\operatorname{relint}$ denotes the relative interior. By \cite[Proposition~5.4]{ReinerRhoades}, the corresponding interior harmonic spaces form the ideal
\[
 \overline{\calH}_P
 =\bigoplus_{m\geq1}y_0^m
 V_{\calZ_m^\circ}.
\]
We define the \emph{interior graded Ehrhart series} by
\[
\overline{\mathrm E}_P(t,q)
=
\sum_{m\geq1}\operatorname{Hilb}(V_{\calZ_m^\circ},q)t^m,
\]
which is the bigraded Hilbert series of $\overline{\calH}_P$.

The analogy with Ehrhart rings led Reiner and Rhoades to formulate the
following conjecture. We use the shift convention $[M(a,b)]_{(m,j)}=M_{(m+a,j+b)}$ and we denote the \emph{canonical module} of $\calH_P$ by $\omega_{\calH_P}$.
For the basic theory of canonical modules, we refer the reader to~\cite[Chapter~3]{BrunsHerzog}.

\begin{conjecture}[{\cite[Conjecture~5.5]{ReinerRhoades}}]\label{conjecture:ReinerRhoades}
For every lattice polytope $P\subseteq\RR^n$, with $\kk=\RR$:
\begin{enumerate}[\rm (i)]
\item $\calH_P\subseteq\RR[y_0,\yb]$ is Noetherian, or equivalently,
      finitely generated over $\RR$;
\item $\calH_P$ is Cohen--Macaulay;
\item There is a shift $(a,b)\in\ZZ^2$ for which
      $\omega_{\calH_P}\cong\overline{\calH}_P(a,b)$ as
      bigraded $\calH_P$-modules.
\end{enumerate}
\end{conjecture}

Cavey~\cite{Cavey} subsequently constructed lattice triangles with non-finitely generated harmonic algebras, disproving assertion~(i) in general. These examples make it necessary to identify classes of polytopes for which the algebraic approach succeeds.

The classes studied by Reiner and Rhoades include one-dimensional lattice polytopes, antiblocking polytopes, standard simplices, cross-polytopes, cubes, and the order and chain polytopes of finite
posets~\cite{ReinerRhoades}.
% For an antiblocking polytope, the harmonic algebra equals its Ehrhart ring. For an order polytope, the harmonic algebra equals the Ehrhart ring of the associated chain polytope, and the two poset polytopes have the same interior harmonic spaces~\cite[Proposition~5.7 and Theorem~5.8]{ReinerRhoades}.

Two further families have recently been investigated. Crowley and Partida~\cite{CrowleyPartida} proved Conjecture~\ref{conjecture:ReinerRhoades} for unimodular zonotopes.
Libman and Miller~\cite{LibmanMiller} studied hypersimplices and, more generally, hyperplane slices of cubes.
They establish rationality with explicit denominators and graded Ehrhart reciprocity for these slices, and prove that the harmonic
algebra of a hypersimplex is generated in dilation degree one~\cite[Theorems~A and~B]{LibmanMiller}.
They leave Cohen--Macaulayness and the canonical module description for hypersimplices as questions~\cite[Section~6]{LibmanMiller}.

In this paper, we study type~$A$ alcoved polytopes, introduced by Lam and Postnikov~\cite{LamPostnikov},
which include order polytopes and, up to integral affine equivalence, hypersimplices. A \emph{type~$A$ alcoved polytope} is a nonempty bounded polytope $P\subseteq\RR^n$ described by inequalities of the form
\[
 b_{ij}\leq x_j-x_i\leq c_{ij},
\]
where $x_0=0$ and integer bounds $b_{ij},c_{ij}$ are specified for each pair of indices $i<j$ in $\{0,\ldots,n\}$. Our aim is to describe the harmonic algebras of these polytopes.

Our first main result determines every initial algebra of $\calH_P$. For a permutation $\sigma$ of $[n]=\{1,\ldots,n\}$, compress the coordinate fibers of $P$ successively in the order $\sigma(1),\ldots,\sigma(n)$, replacing each nonempty fiber $[L,U]$ by $[0,U-L]$. For more details, see Section~\ref{sec:compression}. Let $Q_\sigma$ denote the resulting compressed polytope, and set
\[
 C_\sigma=\RR_{\geq0}(\{1\}\times Q_\sigma),\qquad
 \Gamma_\sigma=C_\sigma\cap\ZZ^{n+1},\qquad
 A_{Q_\sigma}=\kk[\Gamma_\sigma].
\]
We identify $(m,\alpha)\in\Gamma_\sigma$ with $y_0^m\yb^\alpha$. For a monomial order $\preceq$ and a nonzero polynomial $f\in\kk[y_0,\yb]$, let $\ini_\preceq(f)$ denote the initial monomial of $f$ (without its coefficient). For a vector subspace $W\subseteq\kk[y_0,\yb]$, set
\[
\ini_\preceq W=\Span_{\kk}\{\ini_\preceq(f):f\in W\setminus\{0\}\}.
\]

\begin{theorem}\label{thm:initial}
Let $P\subseteq\RR^n$ be a type~$A$ alcoved polytope. Let $\preceq$ be a monomial order on $\kk[y_0,y_1,\ldots,y_n]$, and let $\sigma$ be the permutation with $y_{\sigma(1)}\preceq\cdots\preceq y_{\sigma(n)}$. Then
\[
 \ini_\preceq\calH_P=A_{Q_\sigma}.
\]
Both $\calH_P$ and $A_{Q_\sigma}$ are finitely generated, normal, and Cohen--Macaulay.
\end{theorem}

We remark that the initial algebra depends only on the order of the
variables and is a normal Cohen--Macaulay affine semigroup ring
for every monomial order.
Since the harmonic algebra and its initial algebra have the same bigraded Hilbert series, it follows that $\mathrm E_P(t,q)$ is rational.

Our second main result realizes the canonical module as an explicit ideal of the harmonic algebra.
We set $J_P=(y_1\cdots y_n)\overline{\calH}_P$.

\begin{theorem}\label{thm:maincanonical}
Let $P\subseteq\RR^n$ be a full-dimensional type~$A$ alcoved polytope. Then $J_P$ is a radical canonical ideal of $\calH_P$, and there is a bigraded isomorphism $$\omega_{\calH_P}\cong\overline{\calH}_P(0,-n).$$

Moreover, for every monomial order $\preceq$ on $\kk[y_0,y_1,\ldots,y_n]$, if $\sigma$ is the permutation of $[n]$ such that $y_{\sigma(1)}\preceq\cdots\preceq y_{\sigma(n)}$, then $$\ini_\preceq J_P=\omega_{A_{Q_\sigma}}.$$
\end{theorem}

%Since  the integral affine invariance discussed in Section~\ref{section:basics of harmonic algebras} and the reduction to the full-dimensional case in Section~\ref{sec:compression},
Since integral affine equivalence induces a bigraded isomorphism of harmonic algebras (see Section~\ref{section:basics of harmonic algebras})
and every type~$A$ alcoved polytope $P$ is integrally affinely equivalent to a full-dimensional type~$A$ alcoved polytope in $\RR^{\dim P}$ (see Section~\ref{sec:compression}),
we may assume that $P$ is full-dimensional.
Therefore, Theorems~\ref{thm:initial} and~\ref{thm:maincanonical} establish all three assertions of the
Reiner--Rhoades conjecture for type~$A$ alcoved polytopes.
More precisely, let $P$ be a $d$-dimensional type~$A$ alcoved
polytope. Its harmonic algebra $\calH_P$ is finitely generated,
normal, and Cohen--Macaulay. Moreover, its canonical module is
given by
\[
 \omega_{\calH_P}\cong\overline{\calH}_P(0,-d)
\]
as a bigraded $\calH_P$-module.
In particular, these results resolve the questions on Cohen--Macaulayness and canonical
modules for hypersimplices raised by Libman and Miller~\cite[Section~6]{LibmanMiller}.

Stanley's formula for the Hilbert series of a multigraded canonical module~\cite[Theorem~6.1(ii)]{StanleyHilbertFunctions} then yields the graded Ehrhart reciprocity formula
\[
q^d\overline{\mathrm E}_P(t,q)
=(-1)^{d+1}\mathrm E_P(t^{-1},q^{-1}).
\]

This paper is organized as follows.
Section~\ref{section:basics of harmonic algebras} recalls the basics of harmonic spaces and harmonic algebras. 
Section~\ref{sec:compression} develops coordinate compression for type~$A$ alcoved polytopes. 
In Section~\ref{sec:Description of Initial Algebras}, we construct harmonic bases using laminar lattice zonotopes and prove Theorem~\ref{thm:initial}. Finally, Section~\ref{sec:canonical} studies the canonical ideal and proves Theorem~\ref{thm:maincanonical}.

\section{Basics of harmonic algebras}\label{section:basics of harmonic algebras}
In this section, we recall the basics of harmonic spaces and harmonic algebras.
Let $\kk$ be a field of characteristic zero.
Let $S=\kk[\xb]$ and $\mathbb D=\kk[\yb]$, where
$\xb=(x_1,\ldots,x_n)$ and $\yb=(y_1, \ldots, y_n)$. The action of $S$ on $\mathbb D$ and on
$\widehat{\mathbb D}=\kk[[\yb]]$ is given by
$f\odot g=f(\partial/\partial y_1,\ldots,\partial/\partial y_n)g$.
We define a pairing $\langle f, g \rangle$ as the constant term of $f \odot g$.
For a finite set $\calZ\subseteq\ZZ^n$, let
\begin{equation*}
    \mathbf I(\calZ) = \{f \in S : f(\zb) = 0 \text{ for all } \zb \in \calZ\}
\end{equation*}
be its vanishing ideal in $S$. Write $\tau(f)$ for the top degree
homogeneous component of a nonzero $f\in S$. The \emph{associated graded ideal} of $\mathbf I(\calZ)$ is
\begin{equation*}
    \operatorname{gr}\mathbf I(\calZ) = (\tau(f) : f \in \mathbf I(\calZ) \setminus \{0\}).
\end{equation*}

The \emph{point-orbit ring} is $R(\calZ)=S/\operatorname{gr}\mathbf I(\calZ)$.
The \emph{harmonic space} of $\calZ$ is
\begin{align*}
    V_{\calZ} &= (\operatorname{gr}\mathbf I(\calZ))^\perp \\
    &= \{g \in \mathbb D : \langle f, g \rangle = 0 \text{ for all } f \in \operatorname{gr}\mathbf I(\calZ)\} \\
    &= \{g \in \mathbb D : f \odot g = 0 \text{ for all } f \in \operatorname{gr}\mathbf I(\calZ)\}.
\end{align*}
We note that $\operatorname{Hilb}(V_{\calZ},q) = \operatorname{Hilb}(R(\calZ),q)$.

For subsets $A, B \subseteq \RR^n$, the \emph{Minkowski sum} of $A$ and $B$ is
\begin{equation*}
    A + B = \{\ab + \bb : \ab \in A, \bb \in B\}.
\end{equation*}

The following theorem describes how harmonic spaces behave under
Minkowski sums and inclusions of finite lattice sets.

\begin{theorem}[{\cite[Theorem~1.5 and Lemma~5.3]{ReinerRhoades}}]\label{lem:least}
Let $\calZ,\calZ'\subseteq\ZZ^n$ be finite lattice sets. Then
\[
 V_{\calZ}V_{\calZ'}\subseteq V_{\calZ+\calZ'}.
\]
Moreover, if $\calZ\subseteq\calZ'$, then
$V_{\calZ}\subseteq V_{\calZ'}$.
\end{theorem}

For a lattice polytope $P \subseteq \RR^n$, let $\calZ_m = mP \cap \ZZ^n$. 
The \emph{harmonic algebra} of $P$ is
\begin{equation*}
    \calH_P=\bigoplus_{m\geq0}y_0^m V_{\calZ_m}.
\end{equation*}
Since $\calZ_m+\calZ_{m'}\subseteq\calZ_{m+m'}$, Theorem~\ref{lem:least}
shows that the harmonic algebra of $P$ is closed under multiplication.

We next recall the inverse system of the vanishing ideal
$\mathbf I(\calZ)$.
Since this ideal need not be homogeneous, we work in the formal
power series ring $\widehat{\mathbb D}=\kk[[\yb]]$.
Its inverse system is the $\kk$-vector space
\[
 \mathbf I(\calZ)^\perp
 =
 \left\{
 F\in\widehat{\mathbb D}:
 f\odot F=0
 \text{ for every }f\in\mathbf I(\calZ)
 \right\}.
\]
For a nonzero formal power series $F \in\widehat{\mathbb D}$, let $\beta(F)$ denote its bottom degree homogeneous component.
In other words, for $F=F_r+F_{r+1}+\cdots$, where each $F_j$ is homogeneous of
degree $j$ and $F_r\ne0$, set $\beta(F)=F_r$.
We also set $\beta(0)=0$.

\begin{lemma}[{\cite[Corollary~4.4]{ReinerRhoades}}]\label{lem:tail bottom degree}
For a finite lattice set $\calZ\subseteq\ZZ^n$, one has
\[
 V_{\calZ}
 =
 \Span_{\kk}
 \left\{
 \beta(F):F\in\mathbf I(\calZ)^\perp
 \right\}.
\]
\end{lemma}

For $\zb=(z_1,\ldots,z_n)\in\ZZ^n$, define
\[
 \exp(\zb\cdot\yb)
 =
 \sum_{r\geq0}\frac{(z_1y_1+\cdots+z_ny_n)^r}{r!}
 \in\widehat{\mathbb D},
\]
where $\zb\cdot\yb$ is the usual inner product.

\begin{lemma}[{\cite[Lemma~4.6]{ReinerRhoades}}]\label{lem:exponential k-basis}
For a finite lattice set $\calZ\subseteq\ZZ^n$, the inverse system
$\mathbf I(\calZ)^\perp$ has the $\kk$-basis
\[
 \left\{\exp(\zb\cdot\yb):\zb\in\calZ\right\}.
\]
\end{lemma}

Combining these results gives
\[
 V_{\calZ}
 =
 \Span_{\kk}
 \left\{
 \beta\left(
 \sum_{\zb\in\calZ}c_{\zb}\exp(\zb\cdot\yb)
 \right):
 c_{\zb}\in\kk
 \right\}.
\]
In particular, $V_{\emptyset}=0$.

Suppose that $P'=AP+\mathbf{b}$ with $A\in\operatorname{GL}_n(\ZZ)$ and $\mathbf{b}\in\ZZ^n$.
For a finite lattice set $\calZ\subseteq\ZZ^n$ and $\mathbf{v}\in\ZZ^n$, the identities
\[
 \exp((A\zb+\mathbf{v})\cdot\yb)
 =\exp(\mathbf{v}\cdot\yb)\exp(\zb\cdot A^{\mathsf T}\yb)
\]
and
\[
 \beta\bigl(\exp(\mathbf{v}\cdot\yb)
 F(A^{\mathsf T}\yb)\bigr)
 =\beta(F)(A^{\mathsf T}\yb)
\]
show that the substitution $h(\yb)\mapsto h(A^{\mathsf T}\yb)$ maps $V_{\calZ}$ isomorphically onto $V_{A\calZ+\mathbf{v}}$. Since $mP'=A(mP)+m\mathbf{b}$, extending this substitution by fixing $y_0$ gives a bigraded algebra isomorphism $\calH_P\cong\calH_{P'}$ carrying $\overline{\calH}_P$ onto $\overline{\calH}_{P'}$.

Since adjoining zero coordinates does not change the harmonic spaces and the direction lattices of the affine hulls are saturated, the same conclusion holds for integral affine equivalences between lattice polytopes in different ambient spaces.

\section{Coordinate compression of alcoved polytopes}\label{sec:compression}

In this section, we introduce coordinate compression and study its behavior on type~$A$ alcoved polytopes. Our main goal is to show that successive compressions of the lattice points of a dilated alcoved polytope agree with the lattice points of the corresponding continuously compressed polytope.

Let $e_1,\ldots,e_n$ denote the standard basis vectors of $\RR^n$. For a finite lattice set $\calZ\subseteq\ZZ^n$ and $i\in[n]$, let $C_i(\calZ)$ denote the
\emph{discrete $i$-compression} of $\calZ$ in the sense of Green and Tao~\cite[Proof of Lemma~2.8]{GreenTao}.
Thus, each coordinate fiber containing $N$ lattice points is replaced by $\{0,\ldots,N-1\}$.

For a polytope $K\subseteq\RR^n$, let $c_i(K)$ denote its \emph{continuous $i$-compression} (coordinate Blaschke shaking)~\cite[Section~3]{FreyerLucas}, which replaces each nonempty fiber $[L,U]$ parallel to $e_i$ by $[0,U-L]$.

For a finite lattice set $\calZ\subseteq\ZZ^n$, a polytope $P\subseteq\RR^n$, and a permutation $\sigma$ of $[n]$, define
\[
D_\sigma(\calZ)
  =C_{\sigma(n)}\cdots C_{\sigma(1)}(\calZ),
\qquad
Q_\sigma
  =c_{\sigma(n)}\cdots c_{\sigma(1)}(P).
\]

To handle the parameter-dependent difference inequalities arising from successive coordinate compressions, we first establish the following lemma using Fourier--Motzkin elimination.

\begin{lemma}\label{lem:elimination}
Let $P(p) = \{x \in \RR^n \mid x_{i_r} - x_{j_r} \leq c_r(p), r \in [N]\}$, where $x_0=0, i_r, j_r \in \{0, \ldots, n\}$, and each $c_r(p)$ is an affine function with integer coefficients. Then every coordinate projection of $P(p)$ admits a description of the same form, and every nonempty bounded coordinate fiber has integral endpoints whenever $p$ and the remaining coordinates are integral.
\end{lemma}

\begin{proof}
Fourier--Motzkin elimination preserves the form of the inequalities: eliminating $x_k$ from $x_i-x_k\leq a(p)$ and $x_k-x_j\leq b(p)$ gives $x_i-x_j\leq a(p)+b(p)$. This proves the first assertion. 

When the parameters and remaining coordinates are integral, each endpoint of a nonempty bounded fiber is the maximum or minimum of finitely many integers, and is therefore integral.
\end{proof}

\begin{proposition}\label{prop:compression}
Let $P \subseteq \RR^n$ be a type~$A$ alcoved polytope, and let $\sigma$ be a permutation of $[n]$. Then $Q_\sigma \subseteq \RR_{\geq 0}^n$ is a bounded rational polytope with $\dim Q_\sigma = \dim P$.
Moreover, for every $m \in \ZZ_{\geq 0}$, $D_\sigma(mP \cap \ZZ^n) = mQ_\sigma \cap \ZZ^n$. At each compression step, every nonempty coordinate fiber over integral coordinates has integral endpoints.
\end{proposition}

\begin{proof}
The case $m=0$ is immediate. We may therefore take $m>0$ and, after relabeling, $\sigma=(1,\ldots,n)$.

After the first $i-1$ compressions, put $p=(m,\alpha_1,\ldots,\alpha_{i-1})$. We claim that the remaining coordinates satisfy difference inequalities $x_j-x_k\leq c_{jk}(p)$, as well as inequalities involving only $p$, with each $c_{jk}$ an integer-coefficient linear function. We prove this by induction on $i$.

For the $i$th compression, the fiber endpoints are
\[
L=\max_j\{x_j-c_{ji}(p)\},
\qquad
U=\min_k\{x_k+c_{ik}(p)\}.
\]
Replacing $[L,U]$ by $[0,U-L]$ amounts to
\[
\alpha_i\geq0,\qquad
x_j-x_k\leq
c_{ji}(p)+c_{ik}(p)-\alpha_i
\]
for every lower--upper pair. The inequalities not involving $x_i$ are retained, so the required form is preserved.

Thus each intermediate compressed set is a rational polytope. By Lemma~\ref{lem:elimination}, every nonempty coordinate fiber over integral coordinates has integral endpoints, and hence $C_i(K\cap\ZZ^n)=c_i(K)\cap\ZZ^n$
for each intermediate polytope $K$. Each continuous compression is a piecewise-affine translation along the fibers, preserving boundedness and dimension. After all $n$ compressions, every coordinate is nonnegative.

Applying this identity at each compression step and using $c_i(mK)=m\,c_i(K)$ gives $D_\sigma(mP\cap\ZZ^n)=mQ_\sigma\cap\ZZ^n$.
\end{proof}

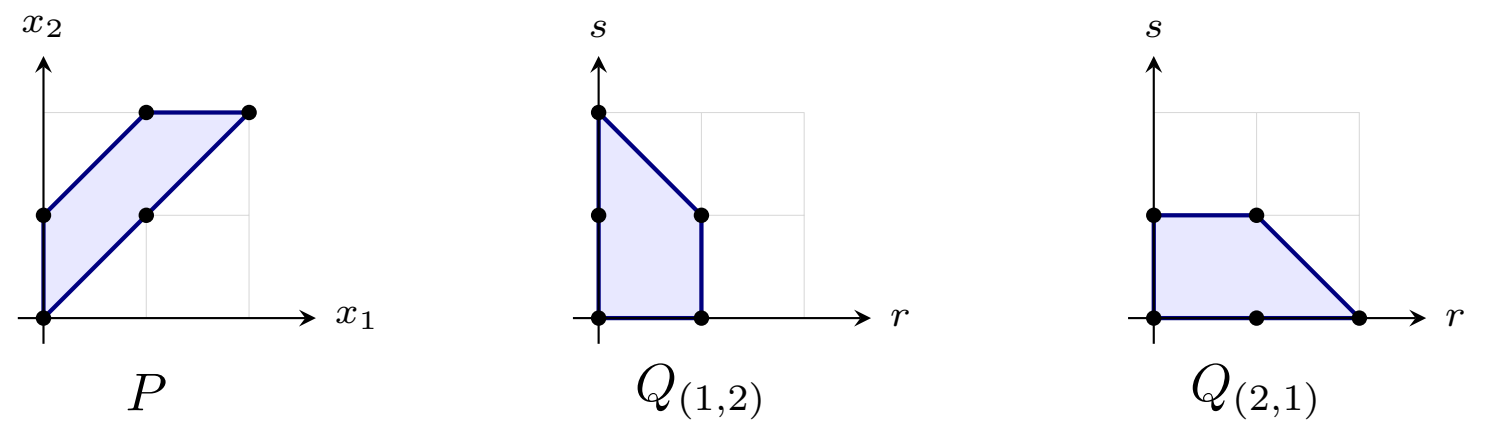
\begin{figure}[htbp]
\centering
\begin{tikzpicture}[x=0.68cm,y=0.68cm,>=stealth,thick]
\foreach \shift/\poly/\points/\xlabel/\ylabel/\captiontext in {
0/{(0,0)--(0,1)--(1,2)--(2,2)}/{0/0,0/1,1/1,1/2,2/2}/x_1/x_2/P,
5.4/{(0,0)--(1,0)--(1,1)--(0,2)}/{0/0,1/0,0/1,1/1,0/2}/r/s/Q_{(1,2)},
10.8/{(0,0)--(2,0)--(1,1)--(0,1)}/{0/0,1/0,2/0,0/1,1/1}/r/s/Q_{(2,1)}}{
\begin{scope}[shift={(\shift,0)}]
 \draw[step=1,gray!30,very thin] (0,0) grid (2,2);
 \filldraw[fill=blue!9,draw=blue!50!black] \poly--cycle;
 \draw[->,thin] (-.25,0)--(2.65,0) node[right] {\scriptsize $\xlabel$};
 \draw[->,thin] (0,-.25)--(0,2.55) node[above] {\scriptsize $\ylabel$};
 \foreach \x/\y in \points \fill (\x,\y) circle (1.45pt);
 \node at (1,-.72) {$\captiontext$};
\end{scope}}
\end{tikzpicture}
\caption{The polytope $P$ and the compressed polytopes $Q_{(1,2)}$ and $Q_{(2,1)}$, together with their lattice points.}
\label{fig:quad}
\end{figure}

\begin{example}\label{ex:compression}
We illustrate successive coordinate compressions and show that the resulting polytope depends on the order of compression. Consider the type~$A$ alcoved polytope
\[
P=\{(x_1,x_2)\in\RR^2:
0\leq x_1\leq x_2\leq 2,\;x_2-x_1\leq 1\}.
\]
Its lattice points are
$P\cap\ZZ^2=\{(0,0),(0,1),(1,1),(1,2),(2,2)\}$.
We compute $Q_{(1,2)}$ and $Q_{(2,1)}$.

First, let $\sigma=(1,2)$. For $m\geq1$, the fiber of $mP$ over $x_2=t$ is
\[
\max\{0,t-m\}\leq x_1\leq t,
\qquad 0\leq t\leq2m.
\]
After compressing the first coordinate, its new coordinate $r$ satisfies $0\leq r\leq\min\{t,m\}$. For fixed $r\in[0,m]$, the second coordinate ranges over $r\leq t\leq2m$. Compressing this coordinate to $s$ gives
\[
mQ_{(1,2)}
=\{(r,s)\in\RR_{\geq0}^2:
r\leq m,\;r+s\leq2m\}.
\]

Next, let $\sigma=(2,1)$. For $m\geq1$, the fiber of $mP$ over $x_1=t$ is
\[
t\leq x_2\leq\min\{t+m,2m\},
\qquad 0\leq t\leq2m.
\]
After compressing the second coordinate, its new coordinate $s$ satisfies $0\leq s\leq\min\{m,2m-t\}$. For fixed $s\in[0,m]$, the first coordinate ranges over $0\leq t\leq2m-s$. Compressing this coordinate to $r$ gives
\[
mQ_{(2,1)}
=\{(r,s)\in\RR_{\geq0}^2:
s\leq m,\;r+s\leq2m\}.
\]

Thus, the compressed polytope depends on the order of compression. When $m=1$, the two compressed polytopes have four lattice points in common: $(0,2)$ belongs only to $Q_{(1,2)}$, whereas $(2,0)$ belongs only to $Q_{(2,1)}$; see Figure~\ref{fig:quad}.
\end{example}

Every type~$A$ alcoved polytope $P$ of dimension $d$ is integrally affinely equivalent to a
full-dimensional type~$A$ alcoved polytope in $\RR^d$. Indeed, the affine hull of $P$ is
defined by those defining inequalities that hold with equality throughout $P$. Each such
equality identifies two coordinates up to an integer constant. These equalities partition
the coordinate indices into equivalence classes. The coordinates in the class containing $x_0$ are fixed. Choosing one representative from each of the other classes identifies the affine lattice of $P$ with $\ZZ^d$. Substituting these representatives into the remaining inequalities again gives difference inequalities.

\section{Description of initial algebras}\label{sec:Description of Initial Algebras}
In this section, we use laminar lattice zonotopes to construct a harmonic basis and determine its initial monomials. For a nonempty set $A\subseteq[n]$, set $\mathbf{1}_A=\sum_{i\in A}e_i$. In particular, $\mathbf{1}=\mathbf{1}_{[n]}=(1,\ldots,1)$. A family of sets is \emph{laminar} if each pair is nested or disjoint. A \emph{laminar lattice zonotope} is a zonotope of the form $Z=\ab+\sum_{A\in\calF}[0,n_A]\mathbf{1}_A$, where $\ab\in\ZZ^n$, $n_A\in\intposit$, and $\calF$ is a laminar family of nonempty subsets of $[n]$. The directions may be linearly dependent. We use the same notation for any finite coordinate index set in place of $[n]$.

\begin{lemma}\label{lem:lift}
Let $P$ be a type~$A$ alcoved polytope with coordinates indexed by $\{1\}\cup J$, where $1\notin J$, and let $\pi$ be the projection that forgets the first coordinate. Then every laminar lattice zonotope $Z'\subseteq\pi(P)$ admits a laminar lattice zonotope lift $Z\subseteq P$ satisfying $\pi(Z)=Z'$.
\end{lemma}

\begin{proof}
Write the segments of $Z'$ as sums of unit segments indexed by a set $E$. For each $e\in E$, let $A_e$ be its support, and write
\[
x_j=a'_j+\sum_{e\in H_j}t_e,
\qquad H_j=\{e\in E:j\in A_e\},
\qquad 0\leq t_e\leq1.
\]
Set $H_0=\emptyset$ and $a'_0=0$. Order the labels by decreasing inclusion of their supports, keeping equal supports consecutive. By laminarity, $H_j\cap H_k$ is an initial segment of both $H_j$ and $H_k$. Indeed, if $e\in H_j\cap H_k$ and $f\in H_j\setminus H_k$, then $A_f\subsetneq A_e$, so $e$ precedes $f$;
the same argument applies to $H_k$.

Write the bounds on the omitted coordinate as
\[
x_j+\lambda_j\leq x_1\leq x_k+\mu_k.
\]
Since $Z'\subseteq\pi(P)$, every lower bound is at most every upper bound for all choices of the parameters. Maximizing their difference gives
\begin{equation}\label{eq:lift-inequality}
a'_j+\lambda_j+|H_j|-(a'_k+\mu_k)
\leq |H_j\cap H_k|.
\end{equation}
Choose $p$ maximizing $a'_j+\lambda_j+|H_j|$, and choose $q$ so that
\[
b=a'_q+\mu_q=\min_k(a'_k+\mu_k).
\]
Put
\[
h=\max\{0,a'_p+\lambda_p+|H_p|-b\}.
\]
By \eqref{eq:lift-inequality}, $h\leq |H_p\cap H_q|$. Let $D$ consist of the first $h$ labels of $H_p\cap H_q$ (so $D=\emptyset$ if $h=0$), and set
\[
x_1=b+\sum_{e\in D}t_e.
\]

Since $D$ is an initial segment of both $H_p$ and $H_q$, we have
\[
|D\cap H_j|=\min\{h,|H_j\cap H_q|\},
\qquad
|D\cap H_k|=\min\{h,|H_p\cap H_k|\}.
\]
The choice of $p,q$, together with \eqref{eq:lift-inequality}, therefore gives
\[
\begin{aligned}
a'_j+\lambda_j+|H_j|-b&\leq |D\cap H_j|,\\
b+h-(a'_k+\mu_k)&\leq |D\cap H_k|.
\end{aligned}
\]
Indeed, the left-hand side of the first inequality is at most both $h$ and $|H_j\cap H_q|$; the second is at most $h$ and, by \eqref{eq:lift-inequality} when $h>0$ (trivially when $h=0$), at most $|H_p\cap H_k|$. As $0\leq t_e\leq1$, it follows that
\[
\begin{aligned}
x_j+\lambda_j-x_1
&\leq a'_j+\lambda_j-b+|H_j\setminus D|\\
&=a'_j+\lambda_j+|H_j|-b-|D\cap H_j|\leq0,\\
x_1-x_k-\mu_k
&\leq b-(a'_k+\mu_k)+|D\setminus H_k|\\
&=b+h-(a'_k+\mu_k)-|D\cap H_k|\leq0.
\end{aligned}
\]
Thus all bounds involving $x_1$ hold, while the remaining inequalities hold because $Z'\subseteq\pi(P)$. Hence the constructed zonotope lies in $P$ and projects onto $Z'$.

Finally, the selected supports form a chain, since they all contain $p$. If a selected support $A_e$ meets an unselected support $A_f$, then $A_f\subseteq A_e$. Otherwise, laminarity gives $A_e\subsetneq A_f$; thus $f\in H_p\cap H_q$ and $f$ precedes $e$, contradicting the choice of $D$. Adjoining $1$ to the selected supports therefore preserves laminarity. Since copies of each support occur consecutively, each original segment is split at most once, into segments of integral lengths. The new translation coordinate $b$ is integral, so the lift is a laminar lattice zonotope.
\end{proof}

\begin{example}\label{ex:split-lift}
We illustrate the lifting construction in Lemma~\ref{lem:lift} using the polytope $P$ from Example~\ref{ex:compression}. Let $\pi$ denote the projection onto the second coordinate. Then $\pi(P)=Z'=[0,2]$.

Splitting $Z'$ into two unit segments, write $x_2=t+u$, where $0\leq t,u\leq1$. The bounds on $x_1$ are
\[
\max\{0,t+u-1\}\leq x_1\leq t+u.
\]
In the notation of the proof of Lemma~\ref{lem:lift}, we have $a_p'+\lambda_p+|H_p|=1$ and $b=0$, so $h=1$. Selecting one of the two unit segments, say the one corresponding to $t$, we set $x_1=t$.

Since $0\leq t\leq t+u\leq2$ and $x_2-x_1=u\leq1$, this gives a laminar lattice zonotope
\[
Z=[0,1](e_1+e_2)+[0,1]e_2 \subseteq P
\]
satisfying $\pi(Z)=Z'$.
\end{example}

For a permutation $\sigma$ let $\max_\sigma A$ be the element of $A$
that occurs last in $\sigma(1),\ldots,\sigma(n)$.

\begin{proposition}\label{prop:witness}
Let $P\subseteq\RR^n$ be a type~$A$ alcoved polytope, let $\sigma$ be a permutation of $[n]$, and set $\calZ=\ZZ^n\cap P$. For every $\alpha\in D_\sigma(\calZ)$ there is a
laminar lattice zonotope $Z_\alpha=\ab+\sum_{A\in\calF}[0,n_A]\mathbf{1}_A
\subseteq P$ such that
\begin{equation}\label{eq:laminar-counts}
               \alpha_i=\sum_{\max_\sigma A=i}n_A\qquad (i\in[n]).
\end{equation}
\end{proposition}
\begin{proof}
After relabeling, we may assume that $\sigma=(1,\ldots,n)$. We proceed by induction on $n$, the case $n=0$ being immediate.

Write $\alpha=(r,\alpha_2,\ldots,\alpha_n)$ and let $\pi$ omit the first coordinate. Consider $E=P\cap(P-re_1)$. Thus $E$ consists of the points $x\in P$ for which $x+re_1$ also belongs to $P$.

By Lemma~\ref{lem:elimination}, the coordinate fibers of $P$ over lattice points have integral endpoints. Hence the first coordinate compression contains $(r,x')$ precisely when $x'\in\pi(E)\cap\ZZ^{n-1}$. Since the remaining compressions act independently on these slices, we obtain
\[
(\alpha_2,\ldots,\alpha_n) \in D_{(2,\ldots,n)}(\pi(E)\cap\ZZ^{n-1}).
\]

By Lemma~\ref{lem:elimination}, $\pi(E)$ is again a type~$A$ alcoved polytope. The induction hypothesis gives a laminar lattice zonotope $Z'\subseteq\pi(E)$ realizing $\alpha_2,\ldots,\alpha_n$. Using the construction in the proof of Lemma~\ref{lem:lift}, we obtain a laminar lattice zonotope $\widetilde Z\subseteq E$. Since this construction only adjoins 1 to some supports, possibly after splitting segments, it preserves the required counts $\alpha_2,\ldots,\alpha_n$.

Now set $Z_\alpha=\widetilde Z+[0,r]e_1$, omitting the added segment if $r=0$. Since both $x$ and $x+re_1$ belong to $P$ for every $x\in\widetilde Z$, convexity gives $Z_\alpha\subseteq P$. The added segment preserves laminarity and contributes exactly $r=\alpha_1$ to the required count. Thus $Z_\alpha$ satisfies all the desired identities.
\end{proof}

\begin{example}\label{ex:witness-construction}
Let $P$ be the polytope in Example~\ref{ex:compression}. We follow the proof of Proposition~\ref{prop:witness} to construct a laminar lattice zonotope for $\sigma=(1,2)$ and $\alpha=(1,1)$.

Since $r=\alpha_1=1$, we first form
\[
 E=P\cap(P-e_1)
   =\{(t,1+t):0\leq t\leq1\}.
\]
Let $\pi$ be the projection onto the second coordinate. Then
\[
 \pi(E)=[1,2],
 \qquad
 T=\pi(E)\cap\ZZ=\{1,2\}.
\]
Compressing $T$ in the second coordinate gives $\{0,1\}$. The one-dimensional induction step therefore realizes $\alpha_2=1$ by the interval
\[
 Z'=[1,2]=\{1+t:0\leq t\leq1\}.
\]

We next lift $Z'$ to $E$ using Lemma~\ref{lem:lift}. Writing $x_2=1+t$, the equation $x_2-x_1=1$ on $E$ forces $x_1=t$. Thus the lift is
\[
 \widetilde Z
 =\{(t,1+t):0\leq t\leq1\}
 =(0,1)+[0,1](e_1+e_2).
\]
In terms of supports, lifting replaces $\{2\}$ by $\{1,2\}$.

Finally, as in the proof of Proposition~\ref{prop:witness}, we adjoin the segment $[0,r]e_1=[0,1]e_1$ and obtain
\[
 Z_\alpha
 =\widetilde Z+[0,1]e_1
 =(0,1)+[0,1](e_1+e_2)+[0,1]e_1.
\]
% Explicitly,
% \[
%  Z_\alpha
%  =\{(t+u,1+t):0\leq t,u\leq1\}.
% \]
% For every point in this zonotope,
% \[
%  0\leq t+u\leq1+t\leq2,
%  \qquad
%  (1+t)-(t+u)=1-u\leq1,
% \]
% so $Z_\alpha\subseteq P$.

The support family is $\calF=\{\{1\},\{1,2\}\}$, with both segment lengths equal to $1$. It is laminar, and
\[
 \max_\sigma\{1\}=1,
 \qquad
 \max_\sigma\{1,2\}=2.
\]
Hence the counts in~\eqref{eq:laminar-counts} are precisely $\alpha_1=\alpha_2=1$.

Similarly, for $\sigma=(1,2)$, we may choose the following laminar lattice zonotopes:
\[
\begin{array}{c|c}
 \alpha & Z_\alpha \\ \hline
 (0,0) & \{(0,0)\} \\[2pt]
 (0,1) & [0,1]e_2 \\[2pt]
 (0,2) & [0,1](e_1+e_2)+[0,1]e_2 \\[2pt]
 (1,0) & (0,1)+[0,1]e_1 \\[2pt]
 (1,1) & (0,1)+[0,1](e_1+e_2)+[0,1]e_1
\end{array}
\]
\end{example}

Now, we construct a $\kk$-basis of $V_{\calZ}$ and $\ini_{\preceq}V_{\calZ}$, where $\calZ=P\cap\ZZ^n$ for a type~$A$ alcoved polytope $P$.
Fix a permutation $\sigma$. For $\alpha\in D_\sigma(\calZ)$, let $Z_{\alpha}=\ab+\sum_{A\in\calF}[0,n_A]\mathbf{1}_A$ be a laminar lattice zonotope satisfying~\eqref{eq:laminar-counts} in Proposition~\ref{prop:witness}.
Let
\begin{equation*}
    h_\alpha=\prod_{A\in\calF}\ell_A(\yb)^{n_A} \text{, where } \ell_A(\yb) = \sum_{i \in A} y_i.
\end{equation*}
% For a monomial order $\preceq$ with $y_{\sigma(1)}\preceq\cdots\preceq y_{\sigma(n)}$, the largest variable in $\ell_A(\yb)$ is $y_{\max_\sigma A}$.
% Therefore, one has $\ini_\preceq(h_\alpha)=\yb^\alpha$.

\begin{proposition}\label{prop:leastbasis}
Let $P\subseteq\RR^n$ be a type~$A$ alcoved polytope and $\sigma$ be a permutation of $[n]$. Set $\calZ=P\cap\ZZ^n$. The set $\{h_\alpha:\alpha\in D_\sigma(\calZ)\}$ forms a $\kk$-basis of $V_{\calZ}$. Moreover, for a monomial order $\preceq$ with $y_{\sigma(1)}\preceq\cdots\preceq y_{\sigma(n)}$, one has $$\ini_\preceq V_{\calZ}=\Span_{\kk}\{\yb^\alpha:\alpha\in D_\sigma(\calZ)\}.$$
\end{proposition}
\begin{proof}
For each $\alpha\in D_\sigma(\calZ)$, let $Z_{\alpha}=\ab+\sum_{A\in\calF}[0,n_A]\mathbf{1}_A\subseteq P$ be a laminar lattice zonotope as in Proposition~\ref{prop:witness} and let $h_\alpha$ and $\ell_A(\yb)$ be as defined immediately before this proposition.
Let $F_{\alpha} = \exp(\ab\cdot\yb)\prod_{A\in\calF}(\exp(\ell_A(\yb))-1)^{n_A}$. We note that $\ell_A(\yb) = \mathbf{1}_A \cdot \yb$.
Expanding $F_\alpha$, we obtain a linear combination of exponentials of the form $\exp((\ab+\sum_{A\in\calF}k_A\mathbf{1}_A)\cdot\yb)$, where $k_A\in\ZZ$ and $0\leq k_A\leq n_A$.
The vector $\ab+\sum_A k_A\mathbf{1}_A$ is a lattice point of
$Z_\alpha$ and hence belongs to $\calZ$, since $Z_\alpha\subseteq P$
by Proposition~\ref{prop:witness}.
Thus, by Lemma~\ref{lem:exponential k-basis}, $F_\alpha\in\mathbf I(\calZ)^\perp$.
The bottom degree homogeneous component of each $\exp(\ell_A(\yb))-1$ is $\ell_A(\yb)$.
Therefore, one has $\beta(F_{\alpha}) = h_\alpha$. By Lemma~\ref{lem:tail bottom degree}, one has $h_\alpha\in V_{\calZ}$.
Let $\preceq$ be a monomial order satisfying $y_{\sigma(1)}\preceq\cdots\preceq y_{\sigma(n)}$.
Then, the largest variable in $\ell_A(\yb)$ is $y_{\max_\sigma A}$. Multiplying initial monomials and using \eqref{eq:laminar-counts} gives $\ini_\preceq(h_\alpha)=\yb^\alpha$. These monomials are distinct, and their number is $|D_\sigma(\calZ)|=|\calZ|=\dim_{\kk}V_{\calZ}$.
Thus, $\{\ini_{\preceq}(h_{\alpha}):\alpha\in D_\sigma(\calZ)\}$ forms a $\kk$-basis of $\ini_{\preceq}V_{\calZ}$.
\end{proof}

\begin{example}\label{ex:quadrilateral-leastbasis}
Retain $P$, $\calZ$, and the laminar lattice zonotopes $Z_\alpha$ from Example~\ref{ex:witness-construction}, with $\sigma=(1,2)$. Fix a monomial order $\preceq$ with $y_1\preceq y_2$. The construction in the proof of Proposition~\ref{prop:leastbasis} gives the following polynomials and their initial monomials:
\[
\begin{array}{c|c|c}
 \alpha & h_\alpha & \ini_\preceq(h_\alpha) \\ \hline
 (0,0) & 1          & 1 \\[2pt]
 (0,1) & y_2        & y_2 \\[2pt]
 (0,2) & y_2 (y_1 + y_2)    & y_2^2 \\[2pt]
 (1,0) & y_1        & y_1 \\[2pt]
 (1,1) & y_1 (y_1 + y_2)    & y_1y_2
\end{array}
\]
For instance, the zonotope $Z_{(1,1)}$ constructed in Example~\ref{ex:witness-construction} gives the exponential sum
\begin{align*}
 F_{(1,1)}
 &= \exp(y_2)
    \bigl(\exp(y_1+y_2)-1\bigr)
    \bigl(\exp(y_1)-1\bigr) \\
 &= \exp(2y_1+2y_2)-\exp(y_1+2y_2)
    -\exp(y_1+y_2)+\exp(y_2).
\end{align*}
Each exponential appearing in this expansion has the form $\exp(\zb\cdot\yb)$ with $\zb\in\calZ$. Hence $F_{(1,1)}\in\mathbf I(\calZ)^\perp$, and
\[
 h_{(1,1)}=\beta(F_{(1,1)})=y_1(y_1+y_2)
\]
belongs to $V_{\calZ}$. The other polynomials are obtained in the same way.

The five polynomials have distinct initial monomials and belong to $V_{\calZ}$. Since $\dim_{\kk}V_{\calZ}=|\calZ|=5$, they form a $\kk$-basis of $V_{\calZ}$. Consequently,
\[
 \ini_\preceq V_{\calZ}
 =\Span_{\kk}\{1,y_1,y_2,y_1y_2,y_2^2\},
\]
whose monomial exponents are precisely the points of $D_{(1,2)}(\calZ)$.
\end{example}

\begin{proof}[Proof of Theorem~\ref{thm:initial}]
Fix $\sigma$ and $\preceq$ satisfying $y_{\sigma(1)}\preceq\cdots\preceq y_{\sigma(n)}$.
Proposition~\ref{prop:leastbasis} applied to $mP$, together with $D_\sigma(mP \cap \ZZ^n) = mQ_\sigma \cap \ZZ^n$ as in
Proposition~\ref{prop:compression}, gives $\ini_\preceq V_{\calZ_m}=\Span_{\kk}\{\yb^\alpha:\alpha\in mQ_\sigma\cap\ZZ^n\}$.
Among the initial monomials of elements of $\calH_P$, those of
$y_0$-degree $m$ are precisely the initial monomials of nonzero
elements of $y_0^m V_{\calZ_m}$.
Since multiplication by $y_0^m$ preserves
the monomial order, we have
\[
 \bigl[\ini_\preceq\calH_P\bigr]_{(m,*)}
 =y_0^m\ini_\preceq V_{\calZ_m}
 =\Span_{\kk}\{y_0^m\yb^\alpha:\alpha\in mQ_\sigma\cap\ZZ^n\}.
\]
Taking the direct sum over $m\geq 0$ gives
$\ini_\preceq\calH_P=\kk[\Gamma_\sigma]$.

The cone $C_\sigma$ is pointed and rational polyhedral by Proposition~\ref{prop:compression}. Gordan's lemma implies that $\Gamma_\sigma$ is finitely generated. Since $\Gamma_\sigma$ contains all the lattice points of $C_\sigma$, it is saturated in the group $\ZZ\Gamma_\sigma$. Its semigroup ring is therefore normal, and it is Cohen--Macaulay by~\cite{Hochster}. Since $\ini_\preceq\calH_P=\kk[\Gamma_\sigma]$ is finitely generated, so is $\calH_P$ by the standard SAGBI basis argument. Its normality and Cohen--Macaulayness follow from~\cite[Theorem~3.16(b)]{BrunsConca2003}.
\end{proof}

\begin{corollary}\label{cor:series}
Let $P\subseteq\RR^n$ be a type~$A$ alcoved polytope. For every permutation $\sigma$ of $[n]$ and every integer $m\geq0$, one has $i_P(m;q)=\sum_{\alpha\in mQ_\sigma\cap\ZZ^n}q^{|\alpha|}$, where $|\alpha|=\alpha_1+\cdots+\alpha_n$. In particular, $\mathrm E_P(t,q)$ is rational.
\end{corollary}
\begin{proof}
    Choose a monomial order $\preceq$ with $y_{\sigma(1)}\preceq\cdots\preceq y_{\sigma(n)}$. The harmonic algebra and its initial algebra have the same bigraded Hilbert series.
    Precisely, one has
    \begin{equation*}
        i_P(m;q)=\operatorname{Hilb}(V_{\calZ_m},q) = \operatorname{Hilb}(\ini_{\preceq}V_{\calZ_m},q) = \sum_{\alpha\in mQ_\sigma\cap\ZZ^n}q^{|\alpha|}.
    \end{equation*}
    Moreover, $\mathrm E_P(t,q)$ is the bigraded Hilbert series of $\kk[\Gamma_\sigma]$ and is therefore rational.
\end{proof}

\section{The canonical ideal}\label{sec:canonical}
Throughout this section, $P\subseteq\RR^n$ is a full-dimensional type~$A$ alcoved polytope, so $n=\dim P$. We use the notation of Section~\ref{sec:compression} for coordinate compression. For $m\geq1$, the relative interior of $mP$ coincides with its interior, so $\calZ_m^\circ=\operatorname{int}(mP)\cap\ZZ^n$.
Recall that
$\overline{\calH}_P=\bigoplus_{m\geq1}y_0^m V_{\calZ_m^\circ}$
and $J_P=(\prod_{i=1}^n y_i)\overline{\calH}_P$.

For $m\geq1$ and $m'\geq0$, we have $\calZ_m^\circ\subseteq\calZ_m$ and $\calZ_m^\circ+\calZ_{m'}\subseteq\calZ_{m+m'}^\circ$.
Theorem~\ref{lem:least} therefore shows that $\overline{\calH}_P$ is an ideal of $\calH_P$.

\begin{lemma}\label{lem:cube-ideal}
$J_P$ is a bihomogeneous ideal of $\calH_P$.
\end{lemma}
\begin{proof}
Fix $m\geq1$ and $\ab\in\calZ_m^\circ$. Let $\varepsilon\in\{0,1\}^n$ and set $a_0=\varepsilon_0=0$. Since $\ab$ is an interior lattice point of $mP$, we have
$mb_{ij}+1\leq a_j-a_i\leq mc_{ij}-1$ for all $0\leq i<j\leq n$.
Since $-1\leq\varepsilon_j-\varepsilon_i\leq1$, we obtain
\[
mb_{ij}\leq(a_j-a_i)+(\varepsilon_j-\varepsilon_i)=(a_j+\varepsilon_j)-(a_i+\varepsilon_i)\leq mc_{ij}.
\]
Thus $\ab+\varepsilon\in mP$, and $\calZ_m^\circ+\{0,1\}^n\subseteq\calZ_m$.

Now consider the exponential sum
\[
F=\prod_{i=1}^n\bigl(\exp(y_i)-1\bigr)=\sum_{\varepsilon\in\{0,1\}^n}(-1)^{n-|\varepsilon|}\exp(\varepsilon\cdot\yb).
\]
The sum on the right belongs to $\mathbf I(\{0,1\}^n)^\perp$ by Lemma~\ref{lem:exponential k-basis}. In the product on the left, the bottom degree homogeneous component of the $i$th factor is $y_i$. Thus $\beta(F)=\prod_{i=1}^n y_i$, and Lemma~\ref{lem:tail bottom degree} gives $\prod_{i=1}^n y_i\in V_{\{0,1\}^n}$. By Theorem~\ref{lem:least},
we have
$\prod_{i=1}^n y_i V_{\calZ_m^\circ}\subseteq V_{\{0,1\}^n+\calZ_m^\circ}\subseteq V_{\calZ_m}$.
Multiplying by $y_0^m$ and summing over $m\geq1$ gives $J_P\subseteq\calH_P$. For $h\in\calH_P$ and $f\in\overline{\calH}_P$, we have $h(f \prod_{i=1}^n y_i )=hf\prod_{i=1}^n y_i\in J_P$ because $hf\in\overline{\calH}_P$. Finally, $J_P$ is bihomogeneous since $\overline{\calH}_P$ is bihomogeneous and $\deg\prod_{i=1}^n y_i=(0,n)$.
\end{proof}

\begin{lemma}\label{lem:open-compression}
For every permutation $\sigma$ of $[n]$ and every $m\geq1$, one has
\[
D_\sigma(\calZ_m^\circ)+\mathbf{1}=\operatorname{int}(mQ_\sigma)\cap\ZZ^n.
\]
\end{lemma}
\begin{proof}
Let $K$ be an intermediate polytope and let $i$ be the next coordinate to be compressed. By Proposition~\ref{prop:compression}, every nonempty $i$th coordinate fiber over integral coordinates has integral endpoints $L$ and $U$. Over a point in the interior of the projection, the lattice points of the open fiber are $L+1,\ldots,U-1$, while those of the compressed fiber are $1,\ldots,U-L-1$. Over a boundary point of the projection, both interior fibers are empty. Thus,
\[
C_i(\operatorname{int}(K)\cap\ZZ^n)+e_i = \operatorname{int}(c_i(K))\cap\ZZ^n.
\]
Since compression in the $i$th coordinate commutes with translations in the other coordinates, iterating this identity gives
\[
D_\sigma(\calZ_m^\circ)+\mathbf{1} = \operatorname{int}(mQ_\sigma)\cap\ZZ^n,
\]
as desired.
\end{proof}

\begin{proposition}\label{prop:initial-canonical}
Let $\preceq$ be a monomial order on $\kk[y_0,\yb]$, and let $\sigma$ be the permutation of $[n]$ with $y_{\sigma(1)}\preceq\cdots\preceq y_{\sigma(n)}$. For every $m\geq1$, one has
\[
\ini_\preceq\bigl([J_P]_{(m,*)}\bigr)=\Span_{\kk}\{y_0^m\yb^\gamma:\gamma\in\operatorname{int}(mQ_\sigma)\cap\ZZ^n\}.
\]
In particular, $\ini_\preceq J_P=\omega_{A_{Q_\sigma}}$.
\end{proposition}
\begin{proof}
Fix $m\geq1$. Let $K_m$ be the closed set defined by $x_0=0$ and
\[
mb_{ij}+1\leq x_j-x_i\leq mc_{ij}-1\qquad0\leq i<j\leq n.
\]
At a lattice point these inequalities are equivalent to $mb_{ij}<x_j-x_i<mc_{ij}$. Hence $K_m\cap\ZZ^n=\calZ_m^\circ$. If $\calZ_m^\circ$ is nonempty, $K_m$ is a type~$A$ alcoved polytope, and Proposition~\ref{prop:leastbasis} applies even if $\dim K_m<n$. Since $y_0^m\prod_{i=1}^n y_i$ is a monomial and $\prod_{i=1}^n y_i=\yb^{\mathbf{1}}$,
\[
\begin{aligned}
\ini_\preceq\bigl([J_P]_{(m,*)}\bigr)
 &=y_0^m\prod_{i=1}^n y_i\,\ini_\preceq V_{\calZ_m^\circ}\\
 &=\Span_{\kk}\{y_0^m\yb^{\alpha+\mathbf{1}}:\alpha\in D_\sigma(\calZ_m^\circ)\}\\
 &=\Span_{\kk}\{y_0^m\yb^\gamma:\gamma\in\operatorname{int}(mQ_\sigma)\cap\ZZ^n\}.
\end{aligned}
\]
The second equality follows from Proposition~\ref{prop:leastbasis} and the third from Lemma~\ref{lem:open-compression}.
If $\calZ_m^\circ$ is empty, both sides vanish by $V_{\emptyset}=0$ and Lemma~\ref{lem:open-compression}.

%By taking the direct sum over $m\geq 0$ on both sides, we have $\ini_\preceq J_P=\omega_{A_{Q_\sigma}}$ by~\cite[Theorem~6.3.5(b)]{BrunsHerzog}.
We now take the direct sum over $m\geq 0$ on both sides. For $m>0$, the point $(m,\gamma)$ belongs to $\operatorname{relint}(C_\sigma)$ exactly when $\gamma\in\operatorname{int}(mQ_\sigma)$. Thus, we have $\ini_\preceq J_P=\omega_{A_{Q_\sigma}}$ by~\cite[Theorem~6.3.5(b)]{BrunsHerzog}.
\end{proof}

\begin{proof}[Proof of Theorem~\ref{thm:maincanonical}]
Fix a monomial order $\preceq$ and let $\sigma$ be as in Proposition~\ref{prop:initial-canonical}. By Theorem~\ref{thm:initial}, $\ini_\preceq\calH_P=A_{Q_\sigma}$, and both rings are finitely generated and Cohen--Macaulay. By Lemma~\ref{lem:cube-ideal} and Proposition~\ref{prop:initial-canonical}, $J_P$ is a bihomogeneous ideal of $\calH_P$ and $\ini_\preceq J_P=\omega_{A_{Q_\sigma}}$. Thus \cite[Theorem~3.19(a)]{BrunsConca2003} shows that $J_P$ is a canonical module of $\calH_P$ after forgetting the grading.

Let $\mathfrak m$ be the homogeneous maximal ideal of $\calH_P$. Since $\mathfrak m$ annihilates $\kk$, it also annihilates $\operatorname{Ext}_{\calH_P}^i(\kk,J_P)$, so these groups may be computed after localization at $\mathfrak m$. The module $(J_P)_{\mathfrak m}$ is a canonical module of the $(n+1)$-dimensional local ring $(\calH_P)_{\mathfrak m}$. Hence \cite[Theorem~3.3.10]{BrunsHerzog} gives $\operatorname{Ext}_{\calH_P}^i(\kk,J_P)=0$ for $i\ne n+1$ and $\dim_{\kk}\operatorname{Ext}_{\calH_P}^{n+1}(\kk,J_P)=1$.

A bigraded free resolution of $\kk$ gives each $\operatorname{Ext}_{\calH_P}^i(\kk,J_P)$ a bigrading. Since $\operatorname{Ext}_{\calH_P}^{n+1}(\kk,J_P)$ is one-dimensional as a $\kk$-vector space, exactly one of its bihomogeneous components is nonzero. Let $(a,b)$ be the bidegree of this component. Any nonzero element of it generates the entire Ext group as a $\calH_P$-module, since $\calH_P$ acts through $\calH_P/\mathfrak m=\kk$.

By our shift convention, $[J_P(a,b)]_{(m,j)}=[J_P]_{(m+a,j+b)}$. Thus a map of bidegree $(a,b)$ from a term of the resolution to $J_P$ has bidegree $(0,0)$ when viewed as a map to $J_P(a,b)$. It follows that the generator of $\operatorname{Ext}_{\calH_P}^{n+1}(\kk,J_P)$ has bidegree $(0,0)$ in $\operatorname{Ext}_{\calH_P}^{n+1}(\kk,J_P(a,b))$. Since the Ext groups vanish for $i\ne n+1$, we obtain isomorphisms of bigraded $\kk$-vector spaces
\[
\operatorname{Ext}_{\calH_P}^i(\kk,J_P(a,b))\cong
\begin{cases}
\kk & \text{if } i=n+1,\\
0 & \text{if } i\ne n+1,
\end{cases}
\]
where $\kk$ on the right is concentrated in bidegree $(0,0)$. Thus $J_P(a,b)$ is a bigraded canonical module by~\cite[Definition~4.1]{BarileBruns}, and $J_P(a,b)\cong\omega_{\calH_P}$ as bigraded $\calH_P$-modules by~\cite[Proposition~4.2(b)]{BarileBruns}.

It remains to determine the shift. Since $\calH_P$ and $A_{Q_\sigma}$ have the same bigraded Hilbert series $\mathrm E_P(t,q)$, Stanley's formula~\cite[Theorem~6.1(ii)]{StanleyHilbertFunctions} gives
\[
\operatorname{Hilb}(\omega_{\calH_P};t,q)=(-1)^{n+1}\mathrm E_P(t^{-1},q^{-1})=\operatorname{Hilb}(\omega_{A_{Q_\sigma}};t,q).
\]
The last series equals $\operatorname{Hilb}(J_P;t,q)$ because $\omega_{A_{Q_\sigma}}=\ini_\preceq J_P$. On the other hand, $J_P(a,b)\cong\omega_{\calH_P}$ gives $\operatorname{Hilb}(\omega_{\calH_P};t,q)=t^{-a}q^{-b}\operatorname{Hilb}(J_P;t,q)$. Hence $a=b=0$. Multiplication by $y_1\cdots y_n$ is injective, has bidegree $(0,n)$, and carries $\overline{\calH}_P$ onto $J_P$. Thus $\overline{\calH}_P(0,-n)\cong J_P\cong\omega_{\calH_P}$ as bigraded $\calH_P$-modules.

For each facet $F$ of $C_\sigma$, the restriction map $A_{Q_\sigma}\to\kk[\Gamma_\sigma\cap F]$ sends every monomial with exponent outside $F$ to zero and fixes the other monomials. Its kernel is prime because the target is a domain. The intersection of these kernels is spanned by the monomials whose exponents lie on no facet, so it equals $\omega_{A_{Q_\sigma}}=\ini_\preceq J_P$. Thus $A_{Q_\sigma}/\ini_\preceq J_P$ is reduced, and~\cite[Theorem~3.16(b)]{BrunsConca2003} implies that $\calH_P/J_P$ is reduced. Hence $J_P$ is radical.
\end{proof}

\section*{Acknowledgements}
The second author was supported by JSPS KAKENHI Grant Number JP25KJ1744. The third author was supported by JST SPRING, Grant Number JPMJSP2138. The authors gratefully acknowledge the assistance of ChatGPT (OpenAI) in exploring proof strategies, locating references, and preparing and revising the manuscript. They take full responsibility for all mathematical statements, proofs, and references.


\begin{thebibliography}{BCRV22}

\bibitem[BB25]{BarileBruns}
M.~Barile and W.~Bruns, \emph{$\ZZ^r$-graded rings and their canonical modules}, in \emph{Commutative Algebra: The Mathematical Legacy of Wolmer V.~Vasconcelos}, De Gruyter, 2025, 389--404.


\bibitem[BC03]{BrunsConca2003}
W.~Bruns and A.~Conca, \emph{Gr\"obner bases and determinantal ideals}, in \emph{Commutative Algebra, Singularities and Computer Algebra}, Kluwer Academic Publishers, 2003, 9--66.

\bibitem[BCRV22]{BrunsConca}
W.~Bruns, A.~Conca, C.~Raicu, and M.~Varbaro, \emph{Determinants, Gr\"obner Bases and Cohomology}, Springer, 2022.

\bibitem[BH98]{BrunsHerzog}
W.~Bruns and J.~Herzog, \emph{Cohen--Macaulay Rings}, revised edition, Cambridge Studies in Advanced Mathematics, vol.~39, Cambridge University Press, 1998.

\bibitem[Cav26]{Cavey}
I.~Cavey, \emph{Graded Ehrhart theory and toric geometry}, Proc. Amer. Math. Soc. \textbf{154} (2026), 1859--1866.

\bibitem[CP26]{CrowleyPartida}
C.~Crowley and E.~Partida, \emph{Graded Ehrhart theory of unimodular zonotopes}, S\'em. Lothar. Combin. \textbf{95B} (2026), Article~41.

\bibitem[FL22]{FreyerLucas}
A.~Freyer and E.~Lucas, \emph{Interpolating between volume and lattice point enumerator with successive minima}, Monatsh. Math. \textbf{198} (2022), 717--740.

\bibitem[GT06]{GreenTao}
B.~Green and T.~Tao, \emph{Compressions, convex geometry and the Freiman--Bilu theorem}, Q. J. Math. \textbf{57} (2006), 495--504.

\bibitem[Hoc72]{Hochster}
M.~Hochster, \emph{Rings of invariants of tori, Cohen--Macaulay rings generated by monomials, and polytopes}, Ann. of Math. (2) \textbf{96} (1972), 318--337.

\bibitem[LP07]{LamPostnikov}
T.~Lam and A.~Postnikov, \emph{Alcoved polytopes, I}, Discrete Comput. Geom. \textbf{38} (2007), 453--478.

\bibitem[LM26]{LibmanMiller}
N.~Libman and W.~Miller, \emph{Graded Ehrhart theory for hypersimplices}, preprint (2026), arXiv:2608.27438.

\bibitem[RR26]{ReinerRhoades}
V.~Reiner and B.~Rhoades, \emph{Harmonics and graded Ehrhart theory}, J. Comb. Algebra, published online (2026), \href{https://doi.org/10.4171/JCA/128}{doi:10.4171/JCA/128}.

\bibitem[Sta78]{StanleyHilbertFunctions}
R.~P.~Stanley, \emph{Hilbert functions of graded algebras}, Adv. Math. \textbf{28} (1978), 57--83.

% \bibitem[Sta86]{StanleyPosetPolytopes}
% R.~P.~Stanley, \emph{Two poset polytopes}, Discrete Comput. Geom. \textbf{1} (1986), 9--23.



\end{thebibliography}
\end{document}